% AMSTeX
\def\ver{Nov. 27, 2001, v.5c}
\documentstyle{amsppt}
\magnification=1200
\hsize=6.5truein
\vsize=8.9truein
\topmatter
\title Monodromy at Infinity and the Weights\\
of Cohomology
\endtitle
\author Alexandru Dimca and Morihiko Saito
\endauthor
\abstract
We show that for a polynomial map, the size of the Jordan blocks for
the eigenvalue 1 of the monodromy at infinity is bounded by the
multiplicity of the reduced divisor at infinity of a good
compactification of a general fiber.
The existence of such Jordan blocks is related to global invariant
cycles of the graded pieces of the weight filtration.
These imply some applications to period integrals.
We also show that such a Jordan block of size greater than 1 for
the graded pieces of the weight filtration is the restriction of a
strictly larger Jordan block for the total cohomology group.
If there are no singularities at infinity, we have a more precise
statement on the monodromy.
\endabstract
\endtopmatter
\tolerance=1000
\baselineskip=12pt
\document
\def\bC{{\Bbb C}}
\def\bD{{\Bbb D}}
\def\bN{{\Bbb N}}
\def\bP{{\Bbb P}}
\def\bQ{{\Bbb Q}}
\def\bZ{{\Bbb Z}}
\def\cD{{\Cal D}}
\def\cH{{\Cal H}}
\def\cM{{\Cal M}}
\def\of{\overline{f}}
\def\oS{\overline{S}}
\def\oX{\overline{X}}
\def\tH{\widetilde{H}}
\def\Aut{\hbox{{\rm Aut}}}
\def\Ker{\hbox{{\rm Ker}}}
\def\Coker{\hbox{{\rm Coker}}}
\def\Im{\hbox{{\rm Im}}}
\def\Coim{\hbox{{\rm Coim}}}
\def\Cone{\hbox{{\rm Cone}}}
\def\Var{\hbox{{\rm Var}}}
\def\MHM{\text{{\rm MHM}}}
\def\MHS{\text{{\rm MHS}}}
\def\DR{\hbox{{\rm DR}}}
\def\supp{\hbox{{\rm supp}}\,}
\def\Gr{\text{{\rm Gr}}}
\def\an{\text{{\rm an}}}
\def\simto{\buildrel\sim\over\to}
\def\SameAuthor{\vrule height3pt depth-2.5pt width1cm}

\noindent
Let
$ X = \bC^{n+1} $,
$ S = \bC $, and
$ f : X \to S $ be a polynomial map.
Set
$ X_{s} = f^{-1}(s) $ for
$ s \in S $.
Then there is a Zariski-open subset
$ U $ of
$ S $ such that the
$ H^{j}(X_{s},\bQ) $ for
$ s \in U $ form a local system on
$ U $.
It is known that the behavior of the local monodromy at infinity of
this local system is rather different from the local monodromy around
the points in
$ S $, see [9], [16], [17], etc.
Among others, it is often observed that the size of the Jordan
blocks for the eigenvalue
$ 1 $ is smaller than the size for the other eigenvalues.
The latter is bounded by
$ j + 1 $ due to (a generalization of) the monodromy theorem,
and this is optimal for the eigenvalues different from
$ 1 $.

For a general
$ s \in U $, let
$ \oX_{s} $ be a good compactification of
$ X_{s} $ such that
$ \oX_{s} $ is smooth and the (reduced) divisor at infinity
$ D_{s} := \oX_{s}\setminus X_{s} $ is a divisor with normal
crossings.
Let
$ m_{s} $ be the maximum of the multiplicity (i.e. the number of local
irreducible components) of
$ D_{s} $.
This is independent of the choice of a general
$ s \in U $.
In this paper, we show

\medskip\noindent
{\bf 0.1.~Theorem.}
{\it The size of the Jordan blocks for the eigenvalue
$ 1 $ of the local monodromy at infinity is bounded by
$ m_{s} $ and also by
$ j $ for
$ j > 0 $.
In particular, this size is
$ 1 $ if
$ X_{s} $ admits a smooth compactification with a smooth divisor
at infinity {\rm (}e.g. if the hypersurface in
$ \bP^{n} $ defined by
the highest degree part of
$ f $ is reduced and smooth{\rm )}.
}

\medskip
More precisely, the size is bounded by the difference
$ m'_{s} $ between the maximal
weight of
$ H^{j}(X_{s},\bQ) $ and
$ j $, see (0.4).
Note that
$ m_{s} \ge m'_{s} $ in general by [7], and we have the strict
inequality, for example if
$ n = 2 $,
$ m_{s} = 2 $ and the dual graph of
$ D_{s} $ has no cycle.

Another interesting fact is that there is a certain condition on
the relation between the Jordan blocks for the eigenvalue
$ 1 $ and the weight filtration
$ W $ of the natural mixed Hodge structure [7] on
$ H^{j}(X_{s},\bQ) $, and such Jordan blocks are closely related to
global invariant cycles of the graded pieces of the weight filtration.
Let
$ G $ be the monodromy group which is the image of the monodromy
representation
$ \pi_{1}(U,s) \to \Aut\,H^{j}(X_{s},\bQ) $.
Note that
$ W $ is stable by the action of
$ G $, because
$ W $ gives the weight filtration on the local system
$ \{ H^{j}(X_{s},\bQ) \}_{s\in U} $ which underlies a variation of
mixed Hodge structures.
(This may be considered to be one of basic examples of geometric
variations of mixed Hodge structures defined on a Zariski-open
subset of
$ \bC $.)
Let
$ T_{\infty} $ denote the monodromy at infinity.
This is an element
of
$ G $, and is defined by choosing a path between
$ s $ and
$ \infty $ in
$ U $.

\medskip\noindent
{\bf 0.2.~Theorem.}
{\it Assume the monodromy at infinity of
$ \Gr_{i}^{W}H^{j}(X_{s},\bQ) $ has a Jordan block of size
$ r \,(> 0) $ for the eigenvalue
$ 1 $.
Then
$ \Gr_{i'}^{W}H^{j}(X_{s},\bQ) $ has a nonzero global invariant
cycle {\rm (}i.e.
$ (\Gr_{i'}^{W}H^{j}(X_{s},\bQ))^{G} \ne 0) $ with
$ i' = i+r+1 \le j+m'_{s} \,(\le\min\{2j,j+m_{s}\}) $, or
$ i' = i > j $ and
$ r = 1 $.
In particular, we have natural isomorphisms
$$
(\Gr_{i}^{W}H^{j}(X_{s},\bQ))^{G} =
(\Gr_{i}^{W}H^{j}(X_{s},\bQ))^{T_{\infty}}
\quad\text{for
$ i - j > m'_{s} - 2 $.}
$$
If
$ i' = i+r+1 $ {\rm (}e.g. if
$ r > 1 ) $, the given Jordan block is the restriction of a strictly
larger Jordan block for the monodromy of
$ H^{j}(X_{s},\bQ) $ to the graded piece
$ \Gr_{i}^{W} $.
}

\medskip
This is a special case of (2.4-5).
If
$ n = 2 $, it implies that the size of the Jordan blocks for the
eigenvalue
$ 1 $ of the monodromy at infinity on
$ \Gr_{i}^{W}H^{j}(X_{s},\bQ) $ is at most
$ 1 $ (compare to the example in [13] mentioned after (0.4) below).
When
$ n = 1 $, (0.2) follows from [9] (because
$ \Gr_{1}^{W}H^{1}(X_{s},\bQ) $ coincides with the cohomology
of a smooth compactification).
The last assertion of (0.2) means that if the restriction of
$ T_{\infty} $ to
$ \Gr_{i}^{W} $ has a Jordan block of size
$ > 1 $ for the eigenvalue
$ 1 $, then there is a strictly larger Jordan block for
$ T_{\infty} $ on
$ H^{j}(X_{s},\bQ) $.
In some special case, this was observed in [10].
Note that the relation between the Jordan blocks of the local monodromies
and the weight filtration is rather complicated in general, and the above
assertion does not follow from the conditions of admissible variation of
mixed Hodge structure.
See (2.9) for an application to period integrals.

If
$ f : \bC^{n+1}\to\bC $ does not have singularities at infinity
(more precisely, if
$ f $ is cohomologically tame [28]), then the situation becomes
quite simple.
We have
$ \tH^{j}(X_{s},\bQ) = 0 $ for
$ s \in U $ and
$ j \ne n $.
Let
$ m(s,\lambda,r) $ denote the number of Jordan blocks of the
local monodromy of
$ \{H^{n}(X_{s},\bQ)\}_{s\in U} $ at
$ s $ with eigenvalue
$ \lambda $ and size
$ r $, and similarly for
$ m'(s,\lambda,r) $ with
$ H^{n}(X_{s},\bQ) $ replaced by
$ \Gr_{n}^{W}H^{n}(X_{s},\bQ) $.
Let
$ r_{i} = \dim\Gr_{i}^{W}H^{n}(X_{s},\bQ) $ for
$ s \in U $.
The following assertion (except for the one about the
monodromy around
$ s \in S $) was obtained in [10], 4.3-5 under an
additional mild assumption.

\medskip\noindent
{\bf 0.3.~Theorem.} {\it
Assume
$ f : \bC^{n+1}\to\bC $ is cohomologically tame [28] with
$ n \ge 1 $.
Then the local systems
$ \{\Gr_{i}^{W}H^{n}(X_{s},\bQ)\}_{s\in U} $ are constant for
$ i \ne n $ and
$$
\aligned
r_{n+1}
&= m(\infty,1,1) \ge
\dim \hbox{\rm IH}^{1}(\bP^{1},R^{n}f_{*}\bQ_{X}|_{U}),
\\
r_{n+r+1}
&= m'(\infty,1,r) = m(\infty,1,r+1)\quad\text{for }r > 0,
\\
m(s,\lambda,r)
&= m'(s,\lambda,r)+ \delta_{\lambda,1}\delta_{r,1}
\sum_{i\ne n}r_{i}\quad\text{for }s\in S\setminus U.
\endaligned
$$
Furthermore, the difference between
$ r_{n+1} $ and
$ \dim \hbox{\rm IH}^{1}(\bP^{1},R^{n}f_{*}\bQ_{X}|_{U}) $
is equal to the length of the direct factor of the perverse sheaf
$ \Gr_{n+1}^{W}{}^{p}R^{n+1}f_{*}\bQ_{X} $ with discrete support.
}

\medskip
Here
$ \hbox{\rm IH}^{1}(\bP^{1},R^{n}f_{*}\bQ_{X}|_{U}) $ is the
intersection cohomology with coefficients in the local system
$ \{H^{n}(X_{s},\bQ)\}_{s\in U} $ on
$ U $, and
$ {}^{p}R^{i}f_{*} = {}^{p}H^{i}f_{*} $ with the notation of [2].
(See also [12], 0.8 for the case
$ n = 1 $.)
Theorem (0.3) means that each Jordan block of size
$ r $ for the eigenvalue
$ 1 $ of the local monodromy of
$ \{\Gr_{n}^{W}H^{n}(X_{s},\bQ)\}_{s\in U} $ at
$ \infty $ has a nontrivial extension with a global section of
$ \{\Gr_{n+r+1}^{W}H^{n}(X_{s},\bQ)\}_{s\in U} $, and gives
a Jordan block of size
$ r + 1 $ of the local monodromy of
$ \{H^{n}(X_{s},\bQ)\}_{s\in U} $.
Otherwise there are no nontrivial extensions between the Jordan
blocks of the graded pieces of the weight filtration.
See (2.6).

In this paper we show that Theorem (0.1) is a special case of
the following assertion on the relation between the local
monodromy and the weights of the cohomology [7]:

\medskip\noindent
{\bf 0.4.~Theorem.}
{\it Let
$ f : X \to S $ be a morphism of complex
algebraic varieties such that
$ \dim X = n+1 $ and
$ S $ is a smooth curve.
Let
$ \oS $ be the smooth compactification of
$ S $.
Let
$ j $ be a positive integer, and
$ r $,
$ r' $ be integers such that
$ r' < r $.
Then the Jordan blocks of the monodromies of
$ H^{j}(X_{s},\bQ) $ around
$ s_{\infty} \in \oS \setminus S $ for the eigenvalue
$ 1 $ have size
$ \le r - r' $ if
$ H^{j+1}(X,\bQ) $ has weights
$ \le j+r $ and
$ H^{j}(X_{s},\bQ) $ for a general
$ s \in S $ has weights in
$ [j+r', j+r] $.
}

\medskip
This means that the restriction on the weights of the
cohomology of the total space and a general fiber implies a certain
restriction on the monodromy at infinity.
The converse of Theorem (0.4) is true in a weak sense if
$ f $ is proper and
$ X $ is smooth so that
$ H^{j}(X_{s},\bQ) $ is pure of weight
$ j $ for a general
$ s $.
See Remark (ii) after (2.3).

For the proof of Theorems, we use the fact that the action of the
nilpotent part
$ N $ of the monodromy on the nearby cycles at infinity (endowed with
the limit mixed Hodge structure) is a morphism of mixed Hodge
structures of type
$ (-1,-1) $ (see [7])
so that the assertion is reduced to the estimate
of the weights of the nearby cycles at infinity.
Then the point is that we can estimate the weights of the {\it limit}
mixed Hodge structure on the nearby cycles in terms of the weights of
the {\it natural} mixed Hodge structure [7,\,III] on the cohomology of
the general fiber and the total space.
We first consider the spectral sequence (1.6.1) which relates
the hypercohomology over
$ S $ of the (perverse) direct image of the constant sheaf by
$ f $ to that of the graded pure pieces of the direct image.
For each graded piece, the weight filtration on its nearby cycles is
given by the monodromy filtration up to a shift.
The weights on the nearby cycles can be estimated by using its higher
direct image by the inclusion to the compactification of
$ S $ (1.4.4), and then using its hypercohomology over
$ S $ (1.5.1).
Thus we can deduce the assertions.
However we do not discuss the
number of the local irreducible components of the deleted fiber which
may be assumed to be a divisor with normal crossings if
$ X $ is smooth.
(The condition on
$ H^{j+1}(X,\bQ) $ can be replaced by that on
$ H^{0}(S,{}^{p}R^{j+1}f_{*}\bQ_{X}) $, see (2.3).)

Theorem (0.1) was proved in the case
$ n = 1 $ by [9], and the estimate by
$ j $ follows from [28], Cor. 10, if
$ f $ is cohomologically tame in the sense of loc. cit.
See also [16], [17].
For
$ n= j = 2 $, there is an example such that the size of a Jordan block
is
$ j $ (see [13], Example (5.3.2):
$ f(x,y,z) = x+y+z+x^2y^2z^2 $).
Note that Theorem (0.1) does not hold for the monodromy around a point of
$ S \setminus U $ (consider e.g.
$ f(x,y) = y^{2} + x^{3} - 3x $),
although the local analogue is true (0.5).

The corresponding local assertion is more or less well-known.
Let
$ f $ be a holomorphic function on a complex manifold
$ X $ (or, more generally, on an analytic space
$ X $ which is a rational homology manifold).
Then we have
$ N^{n+1} = 0 $ on the nearby cycle sheaf
$ \psi_{f}\bQ_{X}[n] $ (by reducing to the normal crossing
case).
See e.g. [26].
Consequently, Jordan blocks of the monodromy on the
$ j $th cohomology of the Milnor fiber at any point of
$ X $ have size
$ \le j+1 $ (by restricting to a generic hyperplane and using
the vanishing of certain relative cohomology [25]).
See also [15].
Restricting to the eigenvalue
$ 1 $ (and to the reduced cohomology),
it is known that the size is bounded by
$ j $ due to J. Steenbrink [33] (in the
isolated singularity case), D. Barlet [1] (for
$ j = n $) and V. Navarro Aznar [27] (in general).

Actually, we can get a slightly better estimate (which is similar
to Theorem (0.1), but is much easier),
when the singularity has certain {\it equisingularity}.

\medskip\noindent
{\bf 0.5.~Proposition.}
{\it Let
$ \{S_{\alpha}\} $ be a Whitney stratification of
$ X $ satisfying Thom's
$ A_{f} $-condition (which exists at least locally by [20]).
Let
$ r = \max\text{\rm codim}\, S_{\alpha} $.
Then the support of the perverse sheaves
$ N^{j}\psi_{f}\bQ_{X}[n] $,
$ N^{j-1}\varphi_{f,1}\bQ_{X}[n] $ have dimension
$ \le n-j $, and
$ N^{j}\psi_{f}\bQ_{X}[n]
= N^{j-1}\varphi_{f,1}\bQ_{X}[n] = 0 $ for
$ j \ge r $.
See (2.8).
}

\medskip
The first assertion on the nearby cycles is equivalent to
$ \dim \supp \Gr_{n-j}^{W}\psi_{f}\bQ[n] \le n-j $
by (1.4.2), and the last assertion to the vanishing of
$ \Gr_{n-j}^{W}\psi_{f}\bQ[n] $ for
$ j \ge r $.
They imply the assertions on the vanishing cycle sheaf
with unipotent monodromy
$ \varphi_{f,1}\bQ_{X}[n] $, because the latter is isomorphic to
$ N\psi_{f,1}\bQ_{X}[n] $ by the sheaf version of
the {\it local invariant cycle theorem} (1.4.5).
This gives also
$ \dim \supp \Gr_{n-j}^{W}\bQ_{X_{0}}[n] \le n-j $, and
$ \Gr_{n-j}^{W}\bQ_{X_{0}}[n] = 0 $ for
$ j \ge r $, where
$ W $ is the weight filtration of the mixed Hodge Module [30].
Note that we can replace
$ r $ by the maximal number of the local irreducible components
of an embedded resolution of
$ f^{-1}(0) $.

In Sect. 1 we review some basic facts from the theory of
mixed Hodge Modules [29], [30], and prove (0.1--6) in Sect. 2.

In this paper, cohomology of a complex algebraic variety means
that of the associated analytic space.

\bigskip\bigskip
\centerline{{\bf 1. Mixed Hodge Modules}}

\bigskip\noindent
{\bf 1.1.}
For a complex algebraic variety
$ X $, let
$ \MHM(X) $ denote the category of mixed Hodge Modules on
$ X $.
See [30, 4.2].
If
$ X $ is smooth, an object
$ \cM $ of
$ \MHM(X) $ consists of
$ ((M,F,W), (K,W), \alpha) $ where
$ (M,F) $ is a filtered
$ \cD_{X} $-Module with the filtration
$ W $,
$ (K,W) $ is a filtered perverse sheaf with rational coefficients on
$ X^{\an} $, and
$ \alpha $ is an isomorphism of perverse sheaves
$ \DR(M) \simeq K\otimes_\bQ \bC $ compatible with
$ W $.
They satisfy several good conditions.
The filtrations
$ F $ and
$ W $ are called respectively the Hodge and weight filtrations.
The category
$ \MHM(X) $ is an abelian categories, and every morphisms is
strictly compatible with the two filtrations
$ (F,W) $ in the strong sense [29].
In general,
$ \MHM(X,\bQ) $ is defined by using closed embeddings of
open subvarieties of
$ X $ into smooth varieties.
See [30], [31].
Then the underlying perverse sheaf
$ K $ of a mixed Hodge Module
$ \cM $ is globally well-defined, and the functor assigning
$ K $ to
$ \cM $ is faithful and exact.

For morphisms
$ f $ of complex algebraic varieties, we can construct canonical
functors
$ f_{*} $,
$ f_{!} $,
$ f^{*} $,
$ f^{!} $ between the derived categories of bounded complexes
of mixed Hodge Modules
$ D^{b}\MHM(X) $ [30, 4.3-4].
We will denote by
$ H^{i} : D^{b}\MHM(X) \to \MHM(X) $ the natural cohomology
functor.

\medskip\noindent
{\bf Remark.}
Let
$ f : X \to Y $ be a morphism of complex algebraic
varieties, and
$ \cM $ a bounded complex of mixed Hodge Module on
$ X $ with a finite decreasing filtration
$ G $.
Then we have a spectral sequence in the category of mixed Hodge Modules
$$
E_{1}^{j,i-j} =
H^{i}f_{*}\Gr_{G}^{j}\cM \Rightarrow
H^{i}f_{*}\cM.
\leqno(1.1.1)
$$
Indeed, the direct image
$ f_{*}\cM $ is represented by a complex of mixed
Hodge Modules endowed with a filtration induced by the filtration
$ G $ on
$ \cM $ by the definition of
direct image [30, 4.3].
In particular,
$ \Gr_{G}^{i} $ commutes with the direct image, and the spectral
sequence follows (see e.g. [7,\,(1.3.1)]).
(We have a similar assertion for the pull-back functor
$ f^{*} $.)

Applying this to the truncation
$ \tau $ on
$ f_{*}\cM $, we get the Leray spectral sequence
$$
E_{2}^{p,q} = H^{p}g_{*}H^{q}f_{*}\cM \Rightarrow
H^{p+q}(gf)_{*}\cM,
\leqno(1.1.2)
$$
as in [7,\,(1.4.8)] for morphisms
$ f : X \to Y $ and
$ g : Y \to Z $.

\medskip\noindent
{\bf 1.2.}
We say that
$ \cM $ is pure of weight
$ r $ if
$ \Gr_{k}^{W}\cM = 0 $ for
$ k \ne r $.
A pure Hodge Module is also called a polarizable Hodge Module.
It admits a decomposition by strict support
$$
\cM = \oplus_{Z}\cM_{Z},
\leqno(1.2.1)
$$
where the direct sum is taken over irreducible closed subvarieties
$ Z $ of
$ X $, and
$ \cM_{Z} $ has support
$ Z $ or
$ \emptyset $, but has no nontrivial subobject or quotient object with
strictly smaller support.
The underlying perverse sheaf
$ K_{Z} $ of
$ \cM_{Z} $ is an intersection complex with local system
coefficients, i.e. there exist a dense open smooth affine subvariety
$ U $ of
$ Z $ with the inclusion
$ j : U \to Z $, and a local system
$ L_{Z} $ on
$ U $ such that
$ K_{Z} $ is the intermediate direct image
$$
j_{!*}(L_{Z}[\dim Z]) := \Im(j_{!}(L_{Z}[\dim Z]) \to
j_{*}(L_{Z}[\dim Z]))
\leqno(1.2.2)
$$
in the sense of [2].
If
$ Z $ is a curve, then
$ K_{Z}[-1] $ is a sheaf in the usual sense, and is isomorphic to
$ j_{*}L_{Z} $ where the direct image is in the usual sense.
In particular, we have
$$
H^{-1}(S,K_{Z}) = H^{0}(U,L_{Z}).
\leqno(1.2.3)
$$

\medskip\noindent
{\bf Remark.}
Let
$ f $ be a proper morphism of complex algebraic varieties, and
$ \cM $ a pure Hodge Modules of weight
$ n $ on
$ X $.
Then the cohomological direct image
$ H^{j}f_{*}\cM $ is pure of weight
$ n + j $.
See [30, 4.5].
A pure Hodge Module is also stable by the intermediate direct image
(1.2.2) due to loc. cit.

\medskip\noindent
{\bf 1.3.}
For
$ X = pt $, we have naturally an equivalence of categories
$$
\MHM(pt) = \MHS^{p},
\leqno(1.3.1)
$$
where the right-hand side denotes the category of graded-polarizable
mixed Hodge structures with rational coefficients [7] (and
$ F^{p} = F_{-p} $).
See [30, (4.2.12)].
So mixed Hodge Modules on
$ pt $ will be identified with graded-polarizable mixed Hodge
structures.

More generally, a mixed Hodge Module
$ ((M,F,W),(K,W);\alpha) $ on a smooth variety
$ X $ such that
$ K[-d] $ is a local system can be naturally identified with an
admissible variation of mixed Hodge structure [22], [34] by replacing
$ K $ with
$ K[-d] $, and
$ W $ with
$ W[-d] $, where
$ d = \dim X $ and
$ (W[-d])_{k} = W_{k+d} $.
See [30, 3.27].
In particular, a polarizable Hodge Module of weight
$ w $ such that the underlying perverse sheaf is a local system
(up to a shift) can be identified with a polarizable variation of Hodge
structure of weight
$ w - d $.

If
$ X $ is smooth, we will denote by
$ \bQ_{X}^{H}[d] \in \MHM(X) $ the pure Hodge Module of weight
$ d $ corresponding to the constant variation of Hodge structure
of type
$ (0,0) $.
In general,
$ \bQ_{X}^{H} $ is defined in the derived category of
bounded complexes of mixed Hodge Modules
$ D^{b}\MHM(X) $.
See [30].
By the direct image under the structure morphism
$ X \to pt $, we get a mixed Hodge structure on the cohomology of
$ X $.
This coincides with Deligne's mixed Hodge structure [7].
See [31].

\medskip\noindent
{\bf 1.4.}
Let
$ g $ be a nonconstant function on
$ X $.
Put
$ Y = g^{-1}(0) $.
Then we have the nearby and vanishing cycle functors
$ \psi_{g} $ and
$ \varphi_{g} $ .
They are exact functors from
$ \MHM(X,\bQ) $ to
$ \MHM(Y,\bQ) $, and correspond to
$ \psi_{g}[-1] $,
$ \varphi_{g}[-1] $ on the underlying perverse sheaves [6].

The semisimple part
$ T_{s} $ of the monodromy acts naturally on
$ \psi_{g}\cM $,
$ \varphi_{g}\cM $, and
$ \Ker(T_{s}-1) $ is denoted by
$ \psi_{g,1}\cM $,
$ \varphi_{g,1}\cM $.
Let
$ N = (2\pi i)^{-1}\log T_{u} $, where
$ T_{u} $ is the unipotent part of the monodromy.
Then
$ N $ gives morphisms of mixed Hodge Modules
$$
N : \psi_{g}\cM \to \psi_{g}\cM(-1),\quad
N : \varphi_{g}\cM \to \varphi_{g}\cM(-1),
\leqno(1.4.1)
$$
where
$ (-1) $ is the Tate twist as in [7] (i.e. the Hodge filtration is
shifted by
$ -1 $ and the weight filtration by
$ 2 $).

The weight filtration is given by the relative monodromy filtration
in the sense of [8, (1.6.13)] and [34].
In particular, if
$ \cM $ is pure of weight
$ k $, then we have isomorphisms
$$
\aligned
N^{i} : \Gr_{k-1+i}^{W}\psi_{g}\cM
&\simto \Gr_{k-1-i}^{W}\psi_{g}\cM(-i),
\\
N^{i} : \Gr_{k+i}^{W}\varphi_{g,1}\cM
&\simto \Gr_{k-i}^{W}\varphi_{g,1}\cM(-i).
\endaligned
\leqno(1.4.2)
$$

Let
$ U = X \setminus Y $ with the inclusion morphisms
$ i : Y \to X $,
$ j : U \to X $.
Since
$ \psi_{g}\cM $ depends only on
$ \cM|_{U} $,
$ \psi_{g}j_{*}\cM $ for
$ \cM \in \MHM(U) $ will be denoted by
$ \psi_{g}\cM $.
By [30, 2.24] we have a canonical isomorphism
$$
i^{*}j_{*}\cM=
\Cone(N : \psi_{g,1}\cM \to \psi_{g,1}\cM(-1)).
\leqno(1.4.3)
$$

Indeed,
$ \Var : \varphi_{g,1}j_{*}\cM \to
\psi_{g,1}\cM(1) $
is an isomorphism (because
$ \Var $ corresponds to the action of
$ t $ on the underlying
$ \cD $-Module, see e.g. [29, 3.4.12]) and
$ \hbox{\rm can} : \psi_{g,1}\cM \to
\varphi_{g,1}j_{*}\cM $ is identified with
$ N $.
(This is related to [34] when
$ X $ is a smooth curve and the local monodromies are unipotent.)

In particular, we have
$$
\Gr_{k+1+i}^{W}H^{0}i^{*}j_{*}\cM
= (P_{N}\Gr_{k-1+i}^{W}\psi_{g,1}\cM)(-1) ,
\leqno(1.4.4)
$$
where the right-hand side denotes the primitive part by the action of
$ N $.
So the dimension of
$ H^{0}i^{*}j_{*}\cM $ coincides with the number of Jordan blocks
for the eigenvalue
$ 1 $ in the case
$ \dim X = 1 $ and
$ g $ is a local coordinate.

\medskip\noindent
{\bf Remarks.}
(i) With the above notation, assume
$ X $ smooth, or more generally,
$ X $ is a rational homology manifold so that
$ \bQ_{X}[\dim X] $ is the intersection complex.
Since
$ Y $ is a locally principal divisor,
$ \bQ_{Y}[n] $ is a perverse sheaf, where
$ n = \dim Y $.
So
$ \bQ_{Y}^{H}[n] $ is a mixed Hodge Module on
$ Y $, and
we have a short exact sequence of mixed Hodge Modules
$$
0 \to \bQ_{Y}^{H}[n] \to
\psi_{g,1}\bQ_{X}^{H}[n+1]
\overset\text{{\rm can}}\to\to
\varphi_{g,1}\bQ_{X}^{H}[n+1] \to 0,
$$
because the cokernel of
$ \text{\rm can} $ corresponds to the maximal quotient object supported on
$ g^{-1}(0) $, and it vanishes in this case.
See [30, (2.4.3)].
Combined with [29, 5.1.7], this implies
$$
\aligned
\bQ_{Y}^{H}[n]
&= \Ker(N : \psi_{g,1}\bQ_{X}^{H}[n+1] \to
\psi_{g,1}\bQ_{X}^{H}[n+1](-1)),
\\
\varphi_{g,1}\bQ_{X}^{H}[n+1]
&=\Coim(N : \psi_{g,1}\bQ_{X}^{H}[n+1] \to
\psi_{g,1}\bQ_{X}^{H}[n+1](-1)).
\endaligned
\leqno(1.4.5)
$$
The first isomorphism may be viewed as the sheaf version of the
{\it local invariant cycle theorem}, and implies the second.
These hold also in the case
$ X $ and
$ g $ are analytic by loc. cit.

(ii) Assume
$ X $ is a smooth curve, and let
$ 0 \in X $ with a local coordinate
$ t $ such that
$ \{0\}=t^{-1}(0) $.
Let
$ \cM $ be an admissible variation of mixed Hodge structure on
$ U $, which is identified with a mixed Hodge Module on
$ U $.
Assume the monodromy around
$ 0 $ is unipotent.
Then the nearby cycles
$ \psi_{t}\cM $ can be defined as in [32], [34]
by extending the Hodge bundles to Deligne's canonical extension [5]
as subbundles, and then restricting them to the fiber at
$ 0 $.
(In general, we have to use the filtration
$ V $ ([21], [26]) indexed by rational numbers as in [29].)

In particular, the set of integers
$ p $ such that
$ \Gr_{F}^{p} \ne 0 $ does not change by
passing to the limit Hodge filtration.
This is used for an estimate of the size of Jordan blocks, see (2.7).

\medskip\noindent
{\bf 1.5.~Lemma.}
{\it Let
$ S $ be a smooth affine curve with the structure morphism
$ a_{S} : S \to pt $, and
$ \oS $ be the smooth compactification of
$ S $ with the inclusions
$ i : \Sigma := \oS\setminus S \to \oS $,
$ j : S \to \oS $.
Let
$ \cM $ be a pure Hodge Module of weight
$ k $ on
$ S $, and let
$ j_{!*}\cM = \Im(j_{!}\cM \to j_{*}
\cM) $, the
intermediate direct image.
Then we have an exact sequence of mixed Hodge structures
$$
0 \to
H^{0}(a_{\oS})_{*}j_{!*}\cM \to
H^{0}(a_{S})_{*}\cM \to
(a_{\Sigma})_{*}H^{0}i^{*}j_{*}\cM \to
H^{1}(a_{\oS})_{*}j_{!*}\cM \to 0,
\leqno(1.5.1)
$$
and an isomorphism
$$
H^{-1}(a_{\oS})_{*}j_{!*}\cM = H^{-1}(a_{S})_{*}\cM.
\leqno(1.5.2)
$$
In particular,
$ H^{-1}(a_{S})_{*}\cM $ is pure of weight
$ k-1 $, and is isomorphic to the dual of
$ H^{1}(a_{\oS})_{*}j_{!*}\cM $ up to a Tate twist.
}

\medskip\noindent
{\it Proof.}
Since
$ S $ is affine, it is enough to show the short exact sequence of mixed
Hodge Modules
$$
0 \to j_{!*}\cM \to j_{*}\cM \to
i_{*}H^{0}i^{*}j_{*}\cM \to 0 .
$$
Indeed,
$ H^{i}(a_{\oS})_{*}j_{!*}\cM $ is pure of weight
$ k + i $ (see Remark after (1.2)), and
$ \cM $ is self-dual up to a Tate twist due to the polarization.

Let
$ K $ be the underlying
perverse sheaf of
$ \cM $.
Replacing
$ \oS $ with an open neighborhood of
$ \oS \setminus S $,
we may assume that
$ \cM $ is a variation of mixed Hodge structure, i.e.
$ K[-1] $ is a local system.
Then the underlying
$ \bQ $-complex of
$ (j_{!*}\cM)[-1] $ is
$ j_{*}(K[-1]) $ where
$ j_{*} $ is the direct image in the usual sense.
So the assertion follows from the distinguished triangle
$$
\to j_{*}(K[-1]) \to \bold{R}j_{*}(K[-1]) \to
R^{1}j_{*}(K[-1]) \to,
$$
which gives the exact sequence of the underlying perverse sheaves
after shifting the complexes by
$ 1 $.

\medskip\noindent
{\bf Remark.}
If every local monodromy of
$ K[-1]|_{U} $ is unipotent, we can prove the assertion by
using [34].

\medskip\noindent
{\bf 1.6.~Lemma.} {\it
With the notation of (1.5), we have spectral sequences of mixed Hodge
structures
$$
\align
E_{1}^{-k,k+i} =
H^{i}(a_{S})_{*}\Gr_{k}^{W}\cM
&\Rightarrow
H^{i}(a_{S})_{*}\cM,
\tag 1.6.1
\\
E_{1}^{-k,k+i} =
(a_{\Sigma})_{*}H^{i}i^{*}j_{*}\Gr_{k}^{W}\cM
&\Rightarrow
(a_{\Sigma})_{*}H^{i}i^{*}j_{*}\cM,
\tag 1.6.2
\endalign
$$
together with a natural morphism of the spectral sequence (1.6.1)
to (1.6.2).
}

\medskip\noindent
{\it Proof.}
The first spectral sequence is clear by (1.1.1), and the argument is
similar for the second.
The morphism of spectral sequences follows from the canonical morphism
$ j_{*}\cM \to i_{*}i^{*}j_{*}\cM $ which is compatible with the
filtration induced by
$ W $ on
$ \cM $.

\medskip\noindent
{\bf Remark.}
The functor
$ i^{*}j_{*} $ calculates the cohomology of the punctured neighborhood
of the points at infinity of
$ S $, see (1.4.3).
For a perverse sheaf
$ K $ (whose restriction to
$ U $ is a local system shifted by
$ 1 $),
$ H^{-1}i^{*}j_{*}K $ and
$ H^{0}i^{*}j_{*}K $ give respectively the local invariant
and coinvariant cycles (i.e. the kernel and cokernel of
$ T_{\infty}-id $).

If the differential
$ d_{r}^{-k,k-1} : E_{r}^{-k-r,k+r-1} \to E_{r}^{-k,k} $ of the
spectral sequence (1.6.2) is nonzero, some element of
$ \Coker(\Gr_{k}^{W}T_{\infty}-id) $ belongs to the image of
$ T_{\infty}-id $, and the corresponding Jordan block for
$ \Gr_{k}^{W}T_{\infty} $ is the restriction of a bigger Jordan
block for
$ T_{\infty} $.

\bigskip\bigskip
\centerline{\bf 2. Proof of Theorems}

\bigskip\noindent
{\bf 2.1.~Proposition.}
{\it Let
$ S $ be a smooth affine curve, and
$ \cM $ a mixed Hodge Module on
$ S $ with the weight filtration
$ W $.
Let
$ U $ be a dense open subvariety of
$ S $ on which
$ \cM $ is a variation of mixed Hodge structure.
Assume
$ H^{0}(a_{S})_{*}\cM $ has weights
$ \le m $.
Then
$ H^{0}(a_{S})_{*}\Gr_{k}^{W}\cM $ has weights
$ \le m $ if
$ \cM|_{U} $ has weights
$ \le m+1 $.
More precisely,
$ H^{-1}(a_{S})_{*}\Gr_{i+1}^{W}\cM \ne 0 $ if
$ \Gr_{i}^{W}H^{0}(a_{S})_{*}\Gr_{k}^{W}\cM \ne 0 $ with
$ i > m $.
}

\medskip\noindent
{\it Proof.}
Consider the spectral sequence (1.6.1).
Since
$ S $ is affine,
$ E_{1}^{p,q} = 0 $ if
$ p+q<-1 $ or
$ p+q>0 $.
By (1.5),
$ E_{1}^{-k,k-1} $ is pure of weight
$ k-1 $, and
$ E_{1}^{-k,k-1} = 0 $ for
$ k-1 > m $ using the decomposition (1.2.1) applied to
$ \Gr_{k}^{W}\cM $, because the direct factor of
$ \Gr_{k}^{W}\cM $ with discrete support does not contribute to
$ H^{-1}(a_{S})_{*}\Gr_{k}^{W}\cM $.
So the assertion is clear by [5].

\medskip\noindent
{\bf 2.2.~Proposition.}
{\it With the above notation, assume
$ \cM $ is a pure Hodge Module of weight
$ k $ on
$ S $.
For
$ s \in \oS \setminus S $, let
$ t $ be a local coordinate at
$ s $.
If
$ H^{0}(a_{S})_{*}\cM $ has weights
$ \le m $, then
$ \psi_{t,1}\cM $ has weights between
$ \min\{k-1, 2k-m\} $ and
$ \max\{k-1, m-2\} $.
Conversely, if the monodromy of
$ \psi_{t,1}\cM $ has a Jordan block of size
$ r $ for the eigenvalue
$ 1 $, then either
$ \Gr_{k+r}^{W}H^{0}(a_{S})_{*}\cM \ne 0 $, or
$ \Gr_{k+r}^{W}H^{1}(a_{\oS})_{*}j_{!*}\cM \ne 0 $ with
$ r = 1 $.
}

\medskip\noindent
{\it Proof.}
We have the symmetry of the weights of nearby cycles by (1.4.2).
So it is enough to estimate the maximal weight for the first assertion,
and it is verified by using the exact sequence (1.5.1) and taking
$ H^{0} $ of (1.4.3).
For the last assertion, the existence of a Jordan block of size
$ r $ corresponds to the nonvanishing of the primitive part of
$ \Gr_{k+r}^{W}(a_{\Sigma})_{*}H^{0}i^{*}j_{*}
\cM $ by (1.4.4).
Then it corresponds by (1.5.1) to that of
$ \Gr_{k+r}^{W}H^{0}(a_{S})_{*}\cM $ or
$ \Gr_{k+r}^{W}H^{1}(a_{\oS})_{*}j_{!*}\cM $.
In the second case, we have
$ r = 1 $ because
$ H^{1}(a_{\oS})_{*}j_{!*}\cM $ is pure of weight
$ k+1 $.
So the assertion follows.

\medskip\noindent
{\bf 2.3.~Proofs of (0.4).}
We may assume
$ S $ connected and then affine (because otherwise
$ \oS = S $).
Let
$ \cM = H^{j+1}f_{*}\bQ_{X}^{H} $ so that
$ {}^{p}R^{j+1}f_{*}\bQ_{X} $ is the underlying perverse sheaf of
$ \cM $ (in particular, its restriction to
$ U $ is
$ (R^{j}f_{*}\bQ_{X}|_{U})[1] $).
Here
$ {}^{p}R^{j+1}f_{*} $ means
$ {}^{p}H^{j+1}\bold{R}f_{*} $.
See [2].

By hypothesis
$ \cM|_{U} $ has weights in
$ [j+r'+1, j+r+1] $.
See (1.3) for the shift of weight.
Furthermore
$ H^{0}(a_{S})_{*}\cM $ has weights
$ \le j+r $ by the exact sequence of mixed Hodge structures
$$
0 \to H^{0}(a_{S})_{*}H^{j+1}f_{*}\bQ_{X}^{H} \to
H^{j+1}(X,\bQ) \to
H^{-1}(a_{S})_{*}H^{j+2}f_{*}\bQ_{X}^{H} \to 0 .
\leqno(2.3.1)
$$
This follows from (1.1.2), because
$ S $ is affine so that
$ E_{2}^{p,q} = 0 $ except for
$ p = -1, 0 $.
So
$ H^{0}(a_{S})_{*}\Gr_{k}^{W}\cM $ has weights
$ \le j+r $ by (2.1).

This implies that the
$ \psi_{t,1}\Gr_{k}^{W}\cM $ for
$ k \in [j + r' + 1, j+r+1] $ (and hence
$ \psi_{t,1}\cM $) have weights in
$ [j-r+2r'+1, j+r] $ using (2.2).
This completes the proof of (0.4).

\medskip\noindent
{\bf Remarks.} (i)
Theorem (0.1) follows from (0.4) by using (2.7) below.

\medskip
(ii) Assume
$ X $ is smooth and
$ f $ is proper.
Then
$ H^{j+1}(X,\bQ) $ has weights
$ \le j+r+1 $ if the Jordan blocks for the eigenvalue
$ 1 $
of the local monodromies of
$ H^{j}(X_{s},\bQ) $ at any points of
$ \oS\setminus S $ have size
$ \le r $.
This follows by the same argument as above using
(1.4.3), (1.5.1--2) and (2.3.1).

\medskip
(iii) Assume
$ f $ is a polynomial map.
Let
$ \cM = H^{j+1}f_{*}\bQ_{X}^{H} $ for
$ j > 0 $.
Then
$ H^{i}(a_{S})_{*}\cM = 0 $ for any
$ i $ by (2.3.1).
Let
$ \bD\cM $ be the dual of
$ \cM $.
Then
$ H^{i}(a_{S})_{!}\bD\cM = 0 $ by the duality, and we get natural
isomorphisms
$$
H^{i}(a_{S})_{*}\bD\cM = (a_{\Sigma})_{*}H^{i}i^{*}j_{*}\bD\cM,
$$
using the distinguished triangle
$ \to j_{!} \to j_{*} \to i_{*}i^{*}j_{*} \to $.
For
$ i = -1 $, this reproves a result of Dimca and N\'emethi [11]:
$$
H_{j}(X_{s},\bQ)^{G} = H_{j}(X_{s},\bQ)^{T_{\infty}}.
$$

\noindent
{\bf 2.4.~Theorem.}
{\it Let
$ f : X \to S $ and
$ \oS $ be as in (0.4).
Assume the monodromy
$ T_{\infty} $ of
$ \Gr_{i}^{W}H^{j}(X_{s},\bQ) $ around
$ s_{\infty} \in \oS \setminus S $ has a Jordan block of size
$ r $ for the eigenvalue
$ 1 $, and
$ H^{j+1}(X,\bQ) $ has weights
$ \le i+r $.
Let
$ m'_{s}+j $ be the maximal weight of
$ H^{j}(X_{s},\bQ) $.
Then
$ \Gr_{i'}^{W}H^{j}(X_{s},\bQ) $ has a nonzero global invariant
cycle {\rm (}i.e.
$ (\Gr_{i'}^{W}H^{j}(X_{s},\bQ))^{G} \ne 0 ) $ with
$ i' = i+r+1 \le j+m'_{s} \,(\le\min\{2j,j+m_{s}\}) $, or
$ i' = i $ and
$ r = 1 $.
In the former case {\rm (}e.g. if
$ r > 1 ) $, the given Jordan block is the restriction of a strictly
larger Jordan block for the monodromy of
$ H^{j}(X_{s},\bQ) $ to the graded piece
$ \Gr_{i}^{W} $.
In the latter case, we have
$ i' = i > j+r' $ if
$ j+r' $ is the minimal weight of
$ H^{j}(X_{s},\bQ) $ and
$ H^{j}(X,\bQ) $ has weights
$ > j+r' $.
}

\medskip\noindent
{\it Proof.}
The argument is similar to (2.3).
The existence of a Jordan block of size
$ r $ corresponds by (2.2) (with
$ k $ replaced by
$ i+1 $) to the nonvanishing of
$$
\Gr_{i+r+1}^{W}H^{0}(a_{S})_{*}\Gr_{i+1}^{W}\cM
\quad\text{or}\quad
\Gr_{i+r+1}^{W}H^{1}(a_{\oS})_{*}j_{!*}\Gr_{i+1}^{W}\cM
$$
with the notation of (1.5) and (2.3).
In the first case, the assertion follows from the last assertion of (2.1)
together with Remark after (1.6) (using (1.2.3)), where
$ i+r+1 \le 2j $ and
$ 2n $ by [7,\,(8.2.4)] (see also (2.7) below).
In the second case, we have
$ r = 1 $ by (2.2), and
$ i > j+r' $ if
$ H^{j}(X,\bQ) $ has weights
$ > j+r' $, because it implies that
$ \{\Gr_{j+r'}^{W}H^{j}(X_{s},\bQ)\}_{s\in U} $ has no nonzero
global section (using (1.6.1)).
So we get the assertion.

\medskip\noindent
{\bf 2.5.~Corollary.}
{\it Let
$ f : X \to S $ and
$ T_{\infty} $ be as above.
Assume
$ H^{j}(X_{s},\bQ) $ and
$ H^{j+1}(X,\bQ) $ have weights
$ \le m $.
Then we have natural isomorphisms
$$
(\Gr_{i}^{W}H^{j}(X_{s},\bQ))^{G} =
(\Gr_{i}^{W}H^{j}(X_{s},\bQ))^{T_{\infty}} \quad
\text{for $ i > m-2 $}
\leqno(2.5.1)
$$
and both are zero if
$ |\oS \setminus S| > 1 $.
}

\medskip\noindent
{\it Proof.}
For
$ i > m - 2 $, the restriction of
$ T_{\infty} $ to the unipotent monodromy part (i.e. the generalized
eigenspace for the eigenvalue 1) of
$ \Gr_{i}^{W}H^{j}(X_{s},\bQ) $ is semisimple by (2.4).
Hence
$ H^{0}i^{*}j_{*}\Gr_{i+1}^{W}\cM $ is pure of weight
$ i+ 2 $, and
$$
\aligned
\oplus_{s_{\infty}\in\Sigma}
\dim (\Gr_{i}^{W}H^{j}(X_{s},\bQ))^{T_{\infty}}
&=
\dim H^{1}(a_{\oS})_{*}j_{!*}\Gr_{i+1}^{W}\cM
\\
&=
\dim (\Gr_{i}^{W}H^{j}(X_{s},\bQ))^{G}
\endaligned
$$
by (1.5.1).
In particular,
$ |\oS \setminus S| = 1 $ if both sides are nonzero.
So the assertion follows.

\medskip\noindent
{\bf Remark.}
If
$ f $ is a polynomial map and a general fiber
$ X_{s} $ admits a smooth compactification such that the divisor
at infinity is smooth, then
$ H^{j}(X_{s},\bQ) $ has weights in
$ [j, j+1] $, and (2.4) implies that
$ 1 $ is not an eigenvalue of the monodromy at infinity of
$ \Gr_{j}^{W}H^{j}(X_{s},\bQ) $ and the size of the Jordan blocks
for the eigenvalue
$ 1 $ of the monodromy of
$ \Gr_{j+1}^{W}H^{j}(X_{s},\bQ) $ is at most
$ 1 $.
In particular, the last assertion holds also for
$ H^{j}(X_{s},\bQ) $.
(Note that it follows also from Theorem (0.1).)

\medskip\noindent
{\bf 2.6.~Proof of (0.3).}
Recall that
$ f : X \to S $ is cohomologically tame [28] if there is an
algebraic compactification
$ \of : \oX \to S $ of
$ f $ such that the support of the (shifted) perverse sheaf
$ \varphi_{\of-c}{}^{p}R^{k}j_{*}(\bQ_{X}[n+1]) $ is contained in
$ X $ for any
$ c \in \bC $ and
$ k \in \bZ $, where
$ j : X \to \oX $ denotes the inclusion and
$ {}^{p}R^{k}j_{*} $ means
$ {}^{p}H^{k}\bold{R}j_{*} $, see [2].
Note that the condition implies that
$ \varphi_{\of-c}{}^{p}R^{k}j_{*}(\bQ_{X}[n+1]) = 0 $ for
$ k \ne 0 $ (because
$ {}^{p}R^{k}j_{*}(\bQ_{X}[n+1]) $ is supported on
$ \oX \setminus X $), and
$ \varphi_{\of-c}{}^{p}R^{0}j_{*}(\bQ_{X}[n+1]) $ has discrete
support.
If furthermore
$ X = \bC^{n+1} $ and
$ S = \bC $, it is easy to show that
$ {}^{p}R^{k}f_{*}\bQ_{X} = 0 $ for
$ k \ne 1, n+1 $.

Let
$ W $ be the weight filtration on the perverse sheaf
$ {}^{p}R^{n+1}f_{*}\bQ_{X} $ coming from the corresponding
mixed Hodge Module.
Note that there is a shift of index by
$ 1 $ between this weight filtration and that on the
cohomology
$ H^{n}(X_{s},\bQ) $, see (1.3).
We first show
$$
\text{
$ \Gr_{r}^{W}{}^{p}R^{n+1}f_{*}\bQ_{X} $ is a constant sheaf if
$ r \ne n+1 $.}
\leqno(2.6.1)
$$

Consider the Leray spectral sequence
$$
E_{2}^{i,k} = {}^{p}R^{i}\of_{*}{}^{p}R^{k}j_{*}(\bQ_{X}[n+1])
\Rightarrow {}^{p}R^{i+k}f_{*}(\bQ_{X}[n+1])
$$
in the category of perverse sheaves on
$ S $.
It underlies a spectral sequence of mixed Hodge Modules and
the functor
$ \cM \to \Gr_{i}^{W}\cM $ is an exact functor of mixed Hodge
Modules.
So it is enough to show that
$ \Gr_{r}^{W}{}^{p}R^{i}\of_{*}{}^{p}R^{k}j_{*}(\bQ_{X}[n+1]) $ are
locally constant sheaves for
$ r \ne n+1 $, because
$ S $ is simply connected.
This is further reduced to the vanishing of the functor
$ \varphi_{t - c} $ applied to these perverse sheaves on
$ S $ for
$ c \in \bC $.
Then using the weight spectral sequence and the commutativity of
the vanishing cycle functor with the direct image under a proper
morphism, the assertion follows from the hypothesis on the
support of the vanishing cycle functor, because
$ \Gr_{r}^{W}{}^{p}R^{0}j_{*}(\bQ_{X}[n+1]) $ is supported on
$ \oX \setminus X $ for
$ r \ne n+1 $.

Now the assertion follows from (2.6.1) together with (1.5-6) as
in the proof of (2.3).
Indeed, letting
$ \cM = H^{0}f_{*}(\bQ_{X}^{H}[n+1]) $ in (1.6.1), we get
$$
E_{1}^{-k,k+i} = 0\quad\text{unless
$ i = -1 $,
$ k > n+1 $ or
$ i = 0 $,
$ k = n+1 $.}
$$
Furthermore
$ E_{1}^{-k,k-1} $ is pure of weight
$ k - 1 $, and
$ E_{\infty}^{-k,k+i} = 0 $ for any
$ i, k $ in (1.6.1) because
$ H^{i}(S,{}^{p}R^{n+1}f_{*}\bQ_{X}) = 0 $ for any
$ i $.

Let
$ \cM' = \Gr_{n+1}^{W}\cM $, and
$ \oS = \bP^{1} $ with the inclusion morphisms
$ i : \{\infty\} \to \oS $,
$ j : S \to \oS $.
Then (1.6) gives commutative diagrams for
$ r \ge 1 $:
$$
\CD
@. H^{-1}(S,\Gr_{n+r+1}^{W}\cM) @>{\sim}>>
H^{-1}i^{*}j_{*}\Gr_{n+r+1}^{W}\cM)
\\
@. @VV{d_{r}}V @VV{d_{r}}V
\\
\Gr_{n+r}^{W}H^{0}(\oS,j_{!*}\cM') @>>>
\Gr_{n+r}^{W}H^{0}(S,\cM') @>>>
\Gr_{n+r}^{W}H^{0}i^{*}j_{*}\cM'
\endCD
$$
where
$ H^{i}(S,\,) $ means
$ H^{i}(a_{S})_{*} $, and the bottom row is
$\Gr_{n+r}^{W} $ of the exact sequence (1.5.1)
with
$ H^{1}(\oS,j_{!*}\cM') = 0 $.
By the above argument the left vertical morphism
$ d_{r} $ is an isomorphism for any
$ r \ge 1 $.
Since
$ H^{0}(\oS,j_{!*}\cM') $ is pure of weight
$ n+1 $, and
$ H^{0}i^{*}j_{*}\cM' $ has weights
$ \ge n + 2 $ by (1.4), we see that the right vertical morphism,
which is induced by
$ d_{r} : E_{r}^{-n-r-1,n+r} \to E_{r}^{-n-1,n+1} $ of (1.6.2), vanishes
for
$ r = 1 $, and is an isomorphism for
$ r > 1 $.
The first vanishing means the splitting of the extension between the
Jordan blocks for the eigenvalue
$ 1 $ of the local monodromy at infinity of
$ \Gr_{n}^{W}H^{n}(X_{s},\bQ) $ and
$ \Gr_{n+1}^{W}H^{n}(X_{s},\bQ) $.
So the assertions on the local monodromy at infinity follows.

The triviality of local extensions at
$ s \in S \setminus U $ follows from the local classification of
perverse sheaves or regular holonomic
$ \cD $-modules ([3], [4]) which implies that locally there are no
nontrivial extensions between intersection complexes with unipotent
local monodromies.

\medskip\noindent
{\bf 2.7.~Generalization of the monodromy theorem.}
Let
$ f : X \to S $ be a morphism of complex algebraic varieties such that
$ \dim S = 1 $.
By Remark (ii) after (1.4), the size of the Jordan blocks of the local
monodromies does not exceed the maximal length of successive numbers
$ p $ such that
$ \Gr_{F}^{p}H^{j}(X_{s},\bC) \ne 0 $,
because the
$ H^{j}(X_{s},\bQ) $ for
$ s \in U $ form an admissible variation of mixed Hodge structure on
a Zariski-open subset
$ U $ of
$ S $
([14], [19], [34], etc.)
The assertion was first shown in [32] when the generic fiber
is proper smooth (see also [23]).
Combined with Remark below, this gives a generalization of the monodromy
theorem (see [18], [24] in the case the generic fiber is
proper smooth).

\medskip\noindent
{\bf Remark.}
Let
$ Y $ be a complex algebraic variety of dimension
$ n $.
Let
$ h^{j,p,q}(Y) = \dim \Gr_{F}^{p}\Gr_{p+q}^{W}H^{j}(Y,\bC) $.
Then by [7,\,(8.2.4)],
$ h^{j,p,q}(Y) = 0 $ except when
$ (p,q) \in [0,j]\times [0,j] $ with
$ j \le n $, or
$ (p,q) \in [j-n,n]\times [j-n,n] $ with
$ j \ge n $.
If
$ Y $ is smooth, we have furthermore
$ h^{j,p,q}(Y) = 0 $ for
$ p + q < j $ by loc. cit.
In particular,
$ H^{j}(Y,\bQ) $ has weights in
$ [j, 2j] $ for
$ j \le n $, and in
$ [j, 2n] $ otherwise.

\medskip\noindent
{\bf 2.8.~Proof of (0.5).}
It is well-known that
$ N^{n+1} = 0 $ on the nearby cycle sheaf
$ \psi_{f}\bQ_{X}[n] $.
See e.g. [26].
This implies
$ N^{n}= 0 $ on the vanishing cycles with unipotent monodromy
$ \varphi_{f,1}\bQ_{X}[n] $ by (1.4.5).
Now we take a Whitney stratification of
$ X $ as in (0.5).
(Here
$ f^{-1}(0) $ is assumed to be a union of strata.)
For each stratum
$ S_{\alpha} $ in
$ f^{-1}(0) $, let
$ X_{\alpha} $ be a transversal space which is a locally closed
complex submanifold of
$ X $.
Applying the above argument to the restriction of
$ f $ to
$ X_{\alpha} $, we get the assertion on the dimension of the support of
$ \Im \,N^{j} $.

\medskip\noindent
{\bf Remarks.}
(i) We can replace
$ r $ in (0.5) by the maximal number of the local irreducible components
of an embedded resolution of
$ f^{-1}(0) $.
In this case, we can prove (0.5), or rather the equivalent assertion
after (0.5), by reducing to the normal crossing case
and then using the calculation of nearby cycle sheaf as in [33]
or [29, (3.6.10)].
(See also [15], [24].)

(ii) Proposition (0.5) implies the assertion that the size of
the Jordan blocks for the eigenvalue
$ 1 $ of the local monodromy on the
$ j $th reduced cohomology of the Milnor fiber is bounded by
$ j $ (see [1], [27], [33]), because a perverse sheaf
$ K $ on an analytic space
$ Y $ satisfies
$ \cH^{i}K = 0 $ for
$ i < -\dim Y $.
This can be verified by induction on
$ \dim Y $ using the transversal space to each stratum with
positive dimension of a Whitney stratification of
$ Y $ and also the long exact sequence associated to
local cohomology.

As an application of Theorem (0.2), we have the following

\medskip\noindent
{\bf 2.9.~Proposition.}
{\it Let
$ t $ be the coordinate of
$ S $.
With the assumption of (0.2), let
$ i' $ be as there.
Then there exists
$ \gamma \in H_{j}(X_{s},\bZ) $ such that, for any algebraic
differential
$ j $-form
$ \omega $ on
$ X $ whose cohomology class in the de Rham cohomology of
the generic fiber has weights
$ \le i' $, the period integral
$ \int_{\gamma_{t}}\omega $ is a (univalent) rational function of
$ t $, where
$ \gamma_{t} $ is a multivalued section of the local system
consisting of the homology groups of general fibers, and is obtained
by the parallel translation of
$ \gamma $ using a local
$ C^{\infty} $ trivialization of the restriction of
$ f $ over
$ U $.
This rational function is nonzero if
$ \omega $ is generic.
}

\medskip\noindent
{\it Proof.}
Let
$ W $ be the dual filtration on
$ H_{j}(X_{s},\bQ) = H^{j}(X_{s},\bQ)^{\vee} $, i.e.
$ W_{-k}H_{j}(X_{s},\bQ) = (H^{j}(X_{s},\bQ)/W_{k-1})^{\vee} $ for
$ k \in \bZ $.
The assumption and (0.2) imply that
$ (\Gr_{-i'}^{W}H_{j}(X_{s},\bQ))^{G} \ne 0 $, because the local system
$ \{\Gr_{i'}^{W}H^{j}(X_{s},\bQ)\} $ is selfdual by the polarization,
and is identified with
$ \{\Gr_{-i'}^{W}H_{j}(X_{s},\bQ)\} $.
Take a nonzero element in
$ (\Gr_{-i'}^{W}H_{j}(X_{s},\bQ))^{G} $, which is represented by
$ \gamma \in W_{-i'}H_{j}(X_{s},\bQ) $.
Then, for an algebraic differential
$ j $-form
$ \omega $ on
$ X $ such that the de Rham cohomology class of its restriction to
the generic fiber of
$ f $ is contained in
$ W_{i'} $, the period integral
$ \int_{\gamma_{t}}\omega $ is univalent,
because the pairing factors through the pairing between
$ \Gr_{-i'}^{W}H_{j}(X_{s},\bQ) $ and
$ \Gr_{i'}^{W}H^{j}(X_{s},\bQ) $.
It is a rational function by regularity,
and is nonzero if
$ \omega $ is generic.
So the assertion follows.

\medskip\noindent
{\bf Remarks.}
(i) Proposition (2.9) does not necessarily imply that
$ \gamma $ is extended to a univalent section of the local system,
because only
$ \Gr_{-i'}^{W}\gamma $ is extended in such a way.
Note that the
$ G $-invariant cycles coincide with the
$ T_{\infty} $-invariant cycles for homology (i.e. for the dual
representation) by Dimca and N\'emethi [11].
However, every invariant cycle of
$ \Gr^{W}H_{j}(X_{s},\bQ) $ does not necessarily come from an
invariant cycle of
$ H_{j}(X_{s},\bQ) $ in general (e.g.
$ f = x^{2}y^{2}z^{2} - x^{2}y^{2} + x^{2} + y^{2}+w^{2} $).

\medskip
(ii) As another application, we have the the following
consequence to the behavior of the period integral at infinity in
general.
For an algebraic differential form
$ \omega $ and
$ \gamma \in H_{j}(X_{s},\bZ) $, consider the asymptotic expansion at
infinity
$$
\int_{\gamma_{t}}\omega \sim
\sum_{\alpha\le\alpha_{0}}\sum_{r=0}^{r(\alpha)}
C(\alpha,r)t^{\alpha}(\log t)^{r},
$$
where
$ \alpha_{0}\in\bQ $,
$ r(\alpha)\in\bN $ and
$ C(\alpha,r)\in\bC $.
Then, by the theory of Nilson class functions in [5],
Theorem (0.1) implies
$$
r(\alpha) \le m'_{s}-1 \,(\le\min\{m_{s}-1,j-1\})\quad
\text{if
$ \alpha \in \bZ $.}
\leqno(2.9.1)
$$
Note that we have only
$ r(\alpha) \le j $ for a general
$ \alpha $ by the monodromy theorem.

\bigskip\bigskip
\centerline{{\bf References}}

\bigskip

\item{[1]}
D.~Barlet, Construction du cup-produit de la fibre de Milnor aux
p\^oles de
$ |f|^{2\lambda} $, Ann. Inst. Fourier, Grenoble 34 (4), (1984),
75--107.

\item{[2]}
A.~Beilinson, J.~Bernstein and P.~Deligne, Faisceaux pervers,
Ast\'erisque, vol. 100, Soc. Math. France, Paris, 1982.

\item{[3]}
L.~Boutet de Monvel, $\cD $-modules holon\^omes r\'eguliers en une
variable, in Math\'ematique et Physiques, Progress in Math.,
Birkh\"auser, vol.~37, (1983), pp.~281--288.

\item{[4]}
J.~Brian\c con and Ph.~Maisonobe, Id\'eaux de germes d'op\'erateurs
diff\'erentiels \`a une variable, Enseign.~Math., 30 (1984), 7--36.

\item{[5]}
P.~Deligne, Equations Diff\'erentielles \`a Points Singuliers
R\'eguliers, Lect. Notes in Math. vol.~163, Springer, Berlin, 1970.

\item{[6]}
\SameAuthor, Le formalisme des cycles \'evanescents, in SGA7 XIII and
XIV, Lect. Notes in Math. vol.~340, Springer, Berlin, 1973, pp. 82--115
and 116--164.

\item{[7]}
\SameAuthor, Th\'eorie de Hodge I, Actes Congr\`es Intern. Math., 1970,
vol. 1, 425-430; II, Publ. Math. IHES, 40 (1971), 5--57; III ibid., 44
(1974), 5--77.

\item{[8]}
\SameAuthor, La conjecture de Weil, Publ. Math. IHES, 52 (1980),
137--252.

\item{[9]}
A.~Dimca, Monodromy at infinity for polynomials in two variables,
Journal of Algebraic Geometry 7 (1998), 771--779.

\item{[10]}
\SameAuthor, Monodromy and Hodge theory of regular functions,
in: New Developments in Singularity Theory (D. Siersma et al. eds.),
Kluwer 2001, 257--278.

\item{[11]}
A.~Dimca and A.~N\'emethi, On the monodromy of complex polynomials,
Duke Math. J. 108 (2001), 199--209.

\item{[12]}
A.~Dimca and M.~Saito, Algebraic Gauss-Manin systems and Brieskorn
modules, Amer. J. Math. 123 (2001), 163--184.

\item{[13]}
A.~Douai, Tr\`es bonnes bases du r\'eseau de Brieskorn d'un polyn\^ome
mod\'er\'e,
Bull. Soc. Math. France 127 (1999), 255--287.

\item{[14]}
F.~El Zein, Th\'eorie de Hodge des cycles \'evanescents.
Ann. Sci. Ecole Norm. Sup. (4) 19 (1986), 107--184.

\item{[15]}
D.~Fried, Monodromy and dynamical systems, Topology 25 (1986),
443--453.

\item{[16]}
R.~Garc{\'\i}a L\'opez and A.~N\'emethi, Hodge numbers attached to a
polynomial map, Ann. Inst. Fourier 49
(1999),1547--1579.

\item{[17]}
\SameAuthor, On the monodromy at infinity of a polynomial map II,
Compos. Math. 115 (1999), 1--20.

\item{[18]}
P.~Griffiths,
Periods of integrals on algebraic manifolds:
Summary of main results and discussion of open problems.
Bull. Amer. Math. Soc. 76 (1970), 228--296.

\item{[19]}
F.~Guill\'en, V.~Navarro Aznar, P.~Pascual-Gainza and F.~Puerta,
Hyperr\'esolutions cubiques et descente cohomologique, Lect. Notes in
Math., Springer, Berlin, vol. 1335, 1988.

\item{[20]}
H.~Hironaka, Stratification and flatness,
in Real and Complex Singularities (Proc. Nordic Summer School, Oslo,
1976) Alphen a/d Rijn: Sijthoff \& Noordhoff 1977, pp. 199--265.

\item{[21]}
M.~Kashiwara, Vanishing cycle sheaves and holonomic systems of
differential equations, Lect. Notes in Math., vol. 1016, Springer, Berlin,
1983, pp. 136--142 .

\item{[22]}
\SameAuthor, A study of variation of mixed Hodge structure,
Publ. RIMS, Kyoto Univ., 22 (1986), 991--1024.

\item{[23]}
N.~Katz,
Nilpotent connections and the monodromy theorem:
Applications of a result of Turrittin,
IHES Publ. Math. 39 (1970), 175--232.

\item{[24]}
A.~Landman,
On the Picard-Lefschetz transformation for algebraic manifolds
acquiring general singularities.
Trans. Amer. Math. Soc. 181 (1973), 89--126.

\item{[25]}
D.T.~L\^e, La monodromie n'a pas de points fixes,
J. Fac. Sci. Univ. Tokyo Sect. IA Math. 22 (1975), 409--427.

\item{[26]}
B.~Malgrange, Polyn\^ome de Bernstein-Sato et cohomologie \'evanescente,
Ast\'erisque, 101--102 (1983), 243--267.

\item{[27]}
V.~Navarro Aznar, Sur la th\'eorie de Hodge-Deligne,
Inv. Math. 90 (1987), 11--76.

\item{[28]}
C.~Sabbah, Hypergeometric periods for a tame polynomial,
C. R. Acad. Sci. Paris, 328 S\'erie I (1999) 603--608.

\item{[29]}
M.~Saito, Modules de Hodge polarisables, Publ. RIMS, Kyoto Univ., 24
(1988), 849--995.

\item{[30]}
\SameAuthor, Mixed Hodge Modules, Publ. RIMS Kyoto Univ. 26 (1990),
221--333.

\item{[31]}
\SameAuthor, Mixed Hodge complexes on algebraic varieties,
Math. Ann. 316 (2000), 283--331.

\item{[32]}
W.~Schmid,
Variation of Hodge structure: the singularities of the period mapping,
Inv. Math. 22 (1973), 211--319.

\item{[33]}
J.~Steenbrink, Mixed Hodge structure on the vanishing cohomology,
in Real and Complex Singularities (Proc. Nordic Summer School, Oslo,
1976) Alphen a/d Rijn: Sijthoff \& Noordhoff 1977, pp. 525--563.

\item{[34]}
J.~Steenbrink and S.~Zucker, Variation of mixed Hodge structure I,
Inv. Math., 80 (1985), 485--542.

\medskip

\noindent
Alexandru Dimca

\noindent
Math\'ematiques pures, Universit\'e Bordeaux I

\noindent
33405 Talence Cedex, FRANCE

\noindent
e-mail: dimca\@math.u-bordeaux.fr

\medskip\noindent
Morihiko Saito

\noindent
RIMS Kyoto University

\noindent
Kyoto 606--8502 JAPAN

\noindent
e-mail: msaito\@kurims.kyoto-u.ac.jp

\medskip\noindent
\ver

\bye